\documentclass[11pt,leqno]{article}

\usepackage{amsfonts, latexsym, amsmath, amssymb, amsthm, hyperref, xcolor}

\usepackage{fullpage}

\usepackage{mdframed, longtable, array, tabu, graphicx, mathrsfs}

\hypersetup{
    colorlinks=true,
    linkcolor=blue,
    filecolor=magenta,      
    urlcolor=cyan,
    citecolor=red,
    pdfpagemode=FullScreen,
    }

\usepackage{cochineal}
\usepackage[cochineal]{newtxmath}

\usepackage[T1]{fontenc}

\newtheorem{theorem}{Theorem}[section]
\newtheorem{lemma}[theorem]{Lemma}

\newtheorem{proposition}[theorem]{Proposition}
\newtheorem{corollary}[theorem]{Corollary}

\theoremstyle{definition}
\newtheorem{definition}[theorem]{Definition}

\theoremstyle{remark}
\newtheorem{remark}[theorem]{Remark}
\newtheorem*{note*}{Note}

\makeatletter

\numberwithin{equation}{section}

\makeatother

\newcommand{\abs}[1]{\left\lvert#1\right\rvert}
\newcommand{\norm}[1]{\left\lVert#1\right\rVert}

\newcommand{\prend}{$\hfill \quad \Box$}

\begin{document}

\title{Dimension-free comparison estimates for suprema of some canonical processes}

\date{}

\author{Shivam Sharma\footnote{Department of Statistics and Data Science, Yale University, New Haven CT 06520, e-mail: \url{shivam.sharma@yale.edu} } }

\maketitle

\begin{abstract}
    \noindent We obtain error approximation bounds between expected suprema of canonical processes that are generated by random vectors with independent coordinates and expected suprema of Gaussian processes. In particular, we obtain a sharper proximity estimate for Rademacher and Gaussian complexities. Our estimates are dimension-free and depend only on the geometric parameters and the numerical complexity of the underlying index set. 
\end{abstract}



\section{Introduction}

\noindent Expected suprema of canonical processes play a very important role in geometric functional analysis and asymptotic convex geometry \cite{AGM-book}, and recently in high-dimensional statistics and machine learning \cite{bartlett2002rademacher}, \cite{RV-recon}, \cite{ALMT_edge}, \cite{7282516}, \cite{oymak2018universality}. On the other hand, even for most prototypical canonical processes (such as Gaussian and Rademacher processes), understanding their expected suprema could be a highly non-trivial endeavour \cite{Tal-reg},\cite{BL-Ber}. Consequently, in many theoretical and practical studies, it is very useful to be able to compare expected suprema of different stochastic processes (without having to estimate them directly). In this paper, we study the class of canonical processes $\{\xi_t\}_{t \in T}$ of the form 
	\[ 
		\xi_t := \langle \xi,t \rangle = \sum_{i=1}^n t_i \xi_i,
	\]
where $\xi$ is a random vector on $\mathbb R^n$ with independent coordinates and $T \subset \mathbb R^n $ is finite. We are interested in the quantity $\mathbb E \sup_{ t\in T} \xi_t $. More precisely, given a canonical process $ \{\xi_t\}_{t \in T}$, we are interested in comparing its expected supremum $ \mathbb E \sup_{t \in T} \xi_t$ with that of the canonical Gaussian process $ \{\langle G,t \rangle \}_{t \in T}$ where $G = (g_1, g_2, \ldots, g_n) $ is a standard Gaussian random vector on $ \mathbb R^n $ (that is, $g_i$ are i.i.d. $\mathcal{N} (0,1)$  random variables for $1 \leq i \leq n$). \\

\noindent In the Gaussian case, the expected behavior of the supremum of a process is very well understood and intimately tied to the geometry of the index set $(T,d)$, where $d$ is the usual Euclidean metric on $ \mathbb R^n $, that is, $d(t,s) := \|t-s\|_2 = (\sum_{i=1}^n |t_i-s_i|^2)^{1/2} $. In fact, the celebrated Fernique-Talagrand majorizing measure theorem \cite{Tal-reg} fully characterizes the expected supremum of a canonical Gaussian process in terms of the geometry of the underlying index set. We refer the reader to \cite[Chapter 2]{Tal-book}, \cite{talagrand2001majorizing}, and \cite{talagrand1996majorizing} for a more general and thorough analysis on suprema of Gaussian processes. \\

\noindent Given that we understand the class of canonical Gaussian processes very well, the main idea is to establish a comparison principle between the expected supremum of a canonical Gaussian process and expected supremum of a canonical process that one wish to understand, thereby leveraging the rich theory already developed about canonical Gaussian processes to gain insights into the behavior of expected suprema of other stochastic processes.\\

\noindent These ideas were used by Talagrand in \cite{Tal-Gibbs} in studying universality in Sherrington-Kirkpatrick (SK) model, a very important model in the theory of spin glasses \cite{Tal-spin}, \cite{CH-univ}. Chatterjee in \cite{chat-error} also used a similar idea for studying error bound in Sudakov-Fernique inequality. Before moving further, let us introduce some notation:
	\begin{align}
		g(T) := \mathbb E \sup_{t \in T} \langle G,t \rangle
	\end{align}
and is called the Gaussian complexity of the set $T$. A random vector on $\mathbb R^n$ with i.i.d. Rademacher (symmetric Bernoulli) coordinates $\epsilon_i$ is denoted by $\mathcal{E} = (\epsilon_1, \epsilon_2, \ldots ,\epsilon_n)$. Let $\mathcal{E}_t := \langle \mathcal{E}, t \rangle = \sum_{i=1}^n t_i \epsilon_i$ \, for $t \in T$. We define
	\begin{align}
		r(T):= \mathbb E \sup_{t \in T} \langle \mathcal{E},t \rangle
	\end{align}
and is called the Rademacher complexity of the set $T$. Here and elsewhere, $C>0$ and $L>0$ denote universal constants, not necessarily the same at each occurrence. \\

\noindent The following comparison ``principle" was proven by Talagrand in \cite{Tal-Gibbs}:
\begin{theorem} [Talagrand] \label{thm:1-1}
	Let $T\subset \mathbb R^n$ be finite, $R_\infty(T) := \sup_{t\in T} \|t\|_\infty$ , and $R_2(T) :=\sup_{t\in T} \|t\|_2$ . Then, there exists an absolute constant C such that
	\[
	\left| r(T) - g(T) \right| \leq C \sqrt{R_2R_\infty} (\log |T|)^{3/4}. 
	\] 
\end{theorem} 

\medskip

\noindent Analogous to the proof of lower bound for Gaussian complexity $g(T)$ in the majorizing measure theorem \cite{talagrand2001majorizing}, the proof of the lower bound for Rademacher complexity $r(T)$ in the Latala-Bednorz theorem\footnote{This was known as \textit{the Bernoulli conjecture} and it took nearly 25 years to prove it \cite[Chapter 6]{Tal-book}.} is based on two fundamental facts, namely, the (sub-gaussian) concentration and the Sudakov minoration for Rademacher processes \cite{BL-Ber}. As in \cite{Tal-Gibbs}, having Theorem \ref{thm:1-1} at our disposal, one can recover the Sudakov minoration for Rademacher processes\footnote{Note that this form is seemingly weaker than the one stated in \cite[Chapter 6]{Tal-book}. However, one can employ this bound ``conditionally'' and derive optimal form of Sudakov minoration for Rademacher processes.} :

\begin{corollary} [Sudakov minoration for Rademacher process]
\noindent There exists a universal constant $C>0$ with the following property. Let $T\subset \mathbb R^n$ be finite and such that $\|u-v\|_2\geq a>0$ for all $u,v\in T$ whenever $u\neq v$, and let $R_2(T)= \sup_{t\in T}\|t\|_2$. Suppose that
\[
\|t\|_\infty \leq Ca^2 / (R_2 \sqrt{\log |T|}), \quad \forall t\in T.
\] Then, we have
\begin{align}
r(T) \geq Ca \sqrt{\log |T|}.
\end{align}    
\end{corollary}

\begin{remark} \label{rem:1-2}
It is well known (see e.g., \cite[Chapters 2 and 6]{Tal-book}) that due to sub-gaussian behavior one has $\abs{g(T)}$, $\abs{r(T)}$ $\leq$ $C R_2 \sqrt{\log |T|}$, and so Theorem \ref{thm:1-1} is of interest 
only when its right-hand side is smaller than this quantity. That is, when $R_{\infty} \sqrt{\log |T|} \leq R_2 $. \\
\end{remark}

\noindent The above Theorem \ref{thm:1-1} together with Remark \ref{rem:1-2} can be intuitively understood as a consequence of the Central Limit Theorem: In general, tails of $\sum_{i=1}^n t_i \epsilon_i$ could be very different from the tails of $\sum_{i=1}^n t_ig_i$. However, when each of the components $t_i$ of $t$ is small compared to $ \norm{t}_2$ (that is, when $R_{\infty} \sqrt{\log |T| } \leq R_2 $), then by the Central Limit Theorem, the random variable $\sum_{i=1}^n t_i\epsilon_i$ will resemble more like a Gaussian random variable. \\

\noindent With the above understanding, one can then look to generalize to other distributions. More precisely, one should expect tighter bounds for $ \abs{ \mathbb E \sup_{t \in T} \langle \xi,t \rangle - g(T)}$ whenever the components  $t_i$ of $t$ are small compared to $\norm{t}_2$ and $R_{\infty} \sqrt{ \log |T| } \leq R_2 $ replaced by a similar condition depending on the distribution of $\xi$. This brings us to our contributions in this paper. First, let us introduce the following notation. For a random vector $\xi = (\xi_1, \xi_2, \ldots , \xi_n)$ on $\mathbb R^n$, define 
	\begin{align}
		c_{\xi}(T) := \mathbb E \sup_{t\in T} \langle \xi,t \rangle = \mathbb E \sup_{ t \in T} \sum_{i=1}^n t_i\xi_i 
	\end{align}
for $t \in T $, where $T \subset \mathbb R^n$ is finite. The quantity $c_\xi(T)$ can be interpreted analogously as measuring geometrical complexity of $T$ with respect to the distribution of the random vector $\xi$.\\

\noindent In this paper, we derive upper estimates for the quantity $\abs{ c_\xi(T) - g(T)}$ under various relaxed assumptions on the random vector $\xi$. 
More precisely:
\begin{itemize}
    \item In Corollary \ref{corr-main}, we derive upper estimates on $\abs{ c_\xi(T) - g(T)}$ under independence and mild moment assumptions on the coordinates of $\xi$, and no further distributional assumption on $\xi$. 
    \item In Theorem \ref{thm-main}, we generalize and sharpen Talagrand's Theorem \ref{thm:1-1} by replacing $\sqrt{R_2 R_{\infty}}$ 
    with a sharper quantity $R_4 (T):= \sup_{t \in T} \norm{t}_4$. 
    \item All our estimates for $\abs{ c_\xi(T) - g(T)}$ depend only on the geometric parameters and numerical complexity (namely, $\log |T|$) of the underlying index set $T$. This is in line with the popular anticipation that ``size'' of a ``nearly'' Gaussian process must be related with the geometry of the underlying metric space $(T,d)$ \cite[Chapter 1]{Tal-book}. Importantly, there is no explicit dependence on the dimension $n$, and hence dimension-free in the title. \\
\end{itemize}

\noindent Our method of proof is inspired by Talagrand's interpolation technique in \cite{Tal-Gibbs} wherein a clever hands-on interpolation is being used to compare $\sum_{i=1}^n t_i\epsilon_i$ 
with $\sum_{i=1}^n t_ig_i$ for $t \in T$. We use a similar idea based on Stein's method by integrating tools from Ornstein-Uhlenbeck semigroup with the tools from theory of spin glasses and Gibbs' measures to construct a more natural and refined interpolation between $c_\xi(T)$ and $g(T)$. Thus, in particular, leading to a sharper estimate in terms of $R_4$ than in Theorem \ref{thm:1-1}in \cite{Tal-Gibbs}.  \\

\noindent Rest of the paper is organized as follows:
\begin{itemize}
\item In Section \ref{sec-2} we provide necessary background material on Gibbs' measures and Ornstein-Uhlenbeck semigroup (Stein's method).
\item In section \ref{sec-3} we present our advertised results (namely Corollary \ref{corr-main} and Theorem \ref{thm-main}). We provide full proofs and discussions on various aspects of the proofs.
\item Section 4: In sub-section \ref{sec-4:picture}, we provide a picture, Figure \ref{error-fig}, facilitating a visual comparison between our Theorem \ref{thm-main} and Talagrand's Theorem \ref{thm:1-1}. Additionally, in \ref{sec-4:examples}, we provide examples witnessing $R_4 \ll   \sqrt{R_2 R_{\infty}}$ in the context of Theorem \ref{thm-main} and Talagrand's Theorem \ref{thm:1-1}. In sub-section \ref{sec-4:univ}, as an application of our results, we discuss universality for \textit{order-m tensors}. In sub-section \ref{sec-4:general}, we discuss generalizations of our results when the index set $T$ is equipped with an arbitrary probability measure $\mu$.  
\end{itemize}

\noindent {\bf Acknowledgements.} The author would like to thank his PhD advisor Petros Valettas for many fruitful discussions and for valuable suggestions regarding the material of this work. The author extends his thanks to the anonymous referee for a careful reading of the manuscript and for his constructive feedback. His valuable comments helped to considerably improve the exposition.


\medskip

\section{Background material} \label{sec-2}

\subsection{Smooth approximation and Gibbs' measures}

Let $\beta>0$ and $T\subset \mathbb R^n $ be finite. Recall the so-called log-partition function \footnote{The motivation for this form of the functions comes from statistical mechanics, where $F_\beta$ is interpreted as the log of a \textit{partition function} and $\beta$ as the \textit{inverse temperature} \cite{Tal-spin}.} $F_\beta$ defined by
\[
 	F_\beta(x) = \frac{1}{\beta} \log \left( \sum_{t\in T} e^{\beta \langle x, t\rangle} \right).
\] 
Partial differentiation with respect to $x_i$ yields
	\begin{align}
		\partial_i F_\beta(x) = \frac{\sum_{t\in T} e^{\beta\langle x, t\rangle} t_i}{\sum_{t\in T} e^{\beta \langle x,t\rangle} }.
	\end{align} 
If we consider the probability measures (called Gibbs' measures in statistical physics; cf. \cite{Bar} for applications in discrete mathematics) $\left \{ \mu_x \mid x \in \mathbb R^n \right \}$ on $T$ assigning mass 
	\begin{align}
		\mu_x(\{t\}) = \frac{e^{\beta\langle x, t\rangle}}{ \sum_{u \in T} e^{\beta\langle x, u\rangle}}, \quad t\in T,
	\end{align} 
then the latter derivative can be expressed as an expectation with respect to this measure, that is,
	\begin{align}
		\partial_i F_\beta(x) = \mathbb E_{\mu_x} [ \ell_i], \quad \ell_i(t) :=\langle t, e_i \rangle=  t_i. 
	\end{align} \\

\noindent We have the following:

\begin{lemma} [Properties of $F_\beta$] \label{lem:prop-Fbeta}
Let $T\subset \mathbb R^n$ be a finite set and let $\beta>0$. Then, we have the following properties:
		\begin{enumerate}
			\item $F_\beta$ is convex.
			\item $F_\beta$ is Lipschitz. For all $x$ we have
				\begin{align*}
					\|\nabla F_\beta(x)\|_2 \leq R_2(T).
				\end{align*}
			\item For all $x\in \mathbb R^n$ one has 
				\begin{align}
					\max_{t \in T} \langle x,t \rangle \leq F_\beta(x) \leq \max_{t \in T} \langle x,t \rangle  +\frac{\log |T|}{\beta}.
				\end{align}
                 In particular, for all $x \in \mathbb R^n $
                 \begin{align*}
                     \abs{F_\beta (x) - \sup_{t \in T} \langle x, t \rangle } \leq \frac{\log \abs{T}}{\beta}.
                 \end{align*}
		\end{enumerate}
		
\end{lemma}

\noindent {\it Proof.} 
\begin{enumerate}	
\item The 2nd order partial derivatives of $F_\beta$ are given by
\[
\partial_{ij}F_\beta(x) = \beta \mathbb E_{\mu_x} [\ell_i \ell_j] - \beta \mathbb E_{\mu_x} [\ell_i] \cdot \mathbb E_{\mu_x} [\ell_j]= \beta [{\rm Cov}(\ell_i , \ell_j)]_{\mu_x}.
\]
\item Let $\|u\|_2=1$. Then, we may write
\[
\langle \nabla F_\beta(x) , u \rangle = \sum_i u_i \partial_i F_\beta(x) = \mathbb E_{\mu_x} [\langle t , u \rangle] \leq \max_{t\in T}\langle t,u \rangle \leq R_2(T).
\] Since $u$ was arbitrary the result follows.
\item Since $\beta > 0$, we may write 
\[
\exp \left( \beta \max_{t\in T} \langle  x, t \rangle \right ) \leq \sum_{t\in T} e^{\beta \langle x,t \rangle} \leq |T| e^{\beta \max_{t\in T} \langle x, t\rangle}.
\]
Applying logarithm and dividing by $\beta$ gives the result. \prend \\
\end{enumerate} 

\noindent Next, we shall obtain some bounds for higher derivatives for $F_\beta$. Similar computations also appeared in \cite[Lemma 2.6]{Tal-Gibbs}.

\begin{lemma}[Bounds for higher derivatives of $F_\beta$] \label{lem:bd-der-Fbeta}
	Let $T \subset \mathbb R^n$ be finite and let $F_\beta$ as above. Then, for all $x\in \mathbb R^n$ we have the following bounds:
	\begin{align}
		|\partial_i^{(3)}F_\beta (x)| \leq 6 \beta^2 \mathbb E_{\mu_x}[|\ell_i|^3], \quad |\partial_i^{(4)}F_\beta (x)| \leq 26 \beta^3 \mathbb E_{\mu_x}[\ell_i^4].
	\end{align} 
	In particular, 
		\begin{align}
			|\partial_i^{(3)}F_\beta (x)| \leq 6 \beta^2 \max_{t\in T}|t_i|^3, \quad |\partial_i^{(4)}F_\beta (x)| \leq 26 \beta^3 \max_{t\in T}|t_i|^4.		
		\end{align}
\end{lemma}

\noindent {\it Proof.} It is straightforward to check that 
	\[
		\partial_i^{(2)} F_\beta(x) = \beta {\rm Var}_{\mu_x}[\ell_i].
	\]

\noindent Similarly, we get 
	\begin{align*}
	  \partial_i^{(3)} F_\beta (x) = \beta^2 \mathbb E_{\mu_x} [\ell_i^3] - 3 \beta^2 \mathbb E_{\mu_x} [\ell_i] \cdot  \mathbb E_{\mu_x} [\ell_i^2] + 2 \beta^2 \left(\mathbb E_{\mu_x} [\ell_i] \right)^3 
								= \beta^2 \mathbb E_{\mu_x} \left[ (\ell_i - \mathbb E_{\mu_x}[\ell_i])^3 \right] , 
	\end{align*}

\noindent and
	\begin{align*}
		\partial_i^{(4)} F_\beta (x) &= \beta^3 \mathbb E_{\mu_x} [\ell_i^4] - 4\beta^3 \mathbb E_{\mu_x} [\ell_i^3] \cdot  \mathbb E_{\mu_x} [\ell_i] - 3 \beta^3 \left(\mathbb E_{\mu_x} [\ell_i^2]\right)^2 +
							 12 \beta^3 \left(\mathbb E_{\mu_x} [\ell_i] \right)^2 \cdot \mathbb E_{\mu_x} [\ell_i^2] - 6 \beta^3 \left( \mathbb E_{\mu_x} [\ell_i] \right)^4 \\
							 & = \beta^3 \mathbb E_{\mu_x} \left[ (\ell_i - \mathbb E_{\mu_x}[\ell_i])^4\right] -3\beta^3 \left( {\rm Var}_{\mu_x}[\ell_i] \right)^2.  
	\end{align*}

\noindent Therefore, triangle inequality and H\"older's inequality yield
	\begin{align*}
        |\partial_i^{(3)}F_\beta (x)| &\leq 6 \beta^2 \mathbb E_{\mu_x}[|\ell_i|^3],
        \end{align*}
\noindent and
        \begin{align*}
        | \partial_i^{(4)} F_\beta (x)| &\le 26\beta^3 \mathbb E_{\mu_x}[\ell_i^4], 
	\end{align*} 
\noindent as desired. 
	
\noindent The ``in particular'' case follows from the simple fact that $\mathbb E_{\mu_x} \left[ \abs{f} \right] \leq \max_{t\in T}|f(t)|$. \prend

\subsection{Stein's method}

\noindent Consider the distance between probability distributions $\mu$ and $\mu_0$ on $\mathbb R^n$ of the form 
\begin{align} \label{prob metric}
 d_{{\cal H}} (\mu, \mu_0) := \sup_{h \in \mathcal{H}} | \mathbb E_{\mu} h - \mathbb E_{\mu_0} h |   ,
\end{align}
where $\mathcal{H}$ is a set of test functions. $\mathcal{H}$ is also called a \textit{separating class} (see \cite[Chapter 3.1 and Appendix C.1]{nourdin2012normal} for a rigorous definition) and such distances are called {\it integral probability metrics}. As described in \cite{Cha-icm}, one of the goals of Stein’s method is to obtain explicit bounds on the distance between a probability distribution $\mu$ of interest and a usually very well-understood approximating distribution $\mu_0$, often called the \textit{target distribution}. A test function $h$ is connected to the distribution of interest through a \textit{Stein Equation}
\begin{align} \label{s-eq}
 h(x) - \mathbb E_{\mu_0} h = \mathcal{T}f(x),   
\end{align}
where $\mathcal{T}$ is a \textit{Stein's operator} for the distribution $\mu_0$ along with an associated \textit{Stein's class} $\mathcal{T(F)}$ of functions such that 
\[
\mathbb E_{\mu_0} [\mathcal{T} f(Z)] = 0 \quad \text{for all} \quad f \in \mathcal{T(F)} \quad \iff  \quad Z \sim \mu_0.
\]
Hence the distance \eqref{prob metric} can be bounded using 
\begin{align*}
    d_{\mathcal{H}} (\mu, \mu_0) \leq \sup_{f \in \mathcal{F}(\mathcal{H})}| \mathbb E_{\mu} \mathcal{T} f(X) |,
\end{align*}
where $\mathcal{F}(\mathcal{H}) = \{ f_h \mid h \in \mathcal{H} \}$ is the set of solutions of the Stein equation \eqref{s-eq} for the test functions $h \in \mathcal{H}$. \\

\noindent For a more in depth and comprehensive understanding of Stein's method, we refer the reader to \cite{Stein1972ABF}, \cite{CGS-book}, \cite{Cha-icm}, and \cite{nourdin2012normal}. Since the basic ingredient for materializing the normal (Gaussian) approximation in the Stein's method is Ornstein-Uhlenbeck (OU) semigroup \cite{BGL-book}, we next recall some of the important definitions and results in this context.

\subsubsection{Ornstein--Uhlenbeck (OU) semigroup}

\begin{definition}[OU semigroup]
      Let $P_t:L^1(\gamma_n) \to L^1(\gamma_n)$ be defined by
\begin{align} 
	(P_tf) (x) := \mathbb E f( e^{-t}x + \sqrt{1-e^{-2t} }G),  
\end{align}
where $ \enspace t\geq 0$, $\enspace G\sim N({\bf 0}, I_n)$, and $\, \gamma_n$ is the standard Gaussian measure on $\mathbb R^n$.
\end{definition}

\noindent The generator $L$ associated with the semigroup $P_t$, that is, $Lf = \lim_{ t \to 0^+} \frac{P_tf -f}{t}$ is given by
\[
Lu(x)= \Delta u (x) - \langle x, \nabla u (x)\rangle 
\]
for $u\in C^2(\mathbb R^n) $. It is known that $P_tf$ solves the following boundary value problem:
\[
\begin{cases}
Lu (x,t) = \partial_t u(x,t) \\ u(x,0)=f(x).			
		\end{cases}
\]

\noindent We also define the following operator:
	\begin{align}
		({\mathscr P} f)(x) := \int_0^\infty \left( P_tf(x)-\mathbb Ef(G) \right) \, dt, \quad x\in \mathbb R^n.
	\end{align}

\noindent Note that the right-hand side converges for Lipschitz maps $f$. Indeed,
	\begin{align*}
	 |P_tf(x) - \mathbb E[f(G)]| &= |\mathbb E [f(e^{-t}x + \sqrt{1-e^{-2t}} G) - f(G)] |\\
	 					&\leq \|f\|_{\rm Lip} \mathbb E \left \|e^{-t}x - \frac{e^{-2t}}{1+\sqrt{1-e^{-2t}}} G \right\| \\
						&\leq e^{-t}\|f\|_{\rm Lip} \left( \|x\| + \mathbb E\|G\| \right),
	\end{align*}
        where $\norm{\cdot}$ is the Euclidean norm. 
        
\hspace{0.2cm}

\noindent Given $f \in L^1(\gamma_n)$ and $x \in \mathbb R^n $,  we may write 
 \begin{align} \label{eq:dtc}
	f(x) - \mathbb E f(G) &= - \int_0^\infty \frac{d}{dt} P_t f(x) \, dt = - \int_0^\infty L P_t f(x) \, = - {\mathscr P} L f (x) = - L {\mathscr P}f (x),
 \end{align}
where the last equality follows from the dominated convergence theorem.

\noindent Also, recall that 
\begin{align} \label{eq:partial-pt}
	\partial_i^{(k)} {\mathscr P} f (x)= \int_0^\infty e^{-tk} P_t (\partial_i^{(k)} f)(x) \, dt.
\end{align}

\noindent For proofs of \eqref{eq:dtc} and \eqref{eq:partial-pt} and for other details on the Ornstein-Uhlenbeck (OU) semigroup, we refer the reader to \cite[Section 2.7.1]{BGL-book}.


\medskip

\section{Main Results} \label{sec-3}

\begin{lemma} \label{lem:s-form}
	Let $\xi = (\xi_1, \ldots, \xi_n)$ be a random vector with independent coordinates, with $\mathbb E\xi_i =0$ and $\mathbb E\xi_i^2=1$. We have the 
	following formulae:
		\begin{enumerate} \label{eq:3-1} 
			\item If $\mathbb E|\xi_i|^3 < \infty$ for all $i\leq n$, then for any $C^3$-smooth function $f: \mathbb R^n\to \mathbb R$ we have
				\begin{align}
					\mathbb E \left[ Lf(\xi)\right] = \sum_{i=1}^n \mathbb E \left[\xi_i \int_0^1 \partial_i^{(3)} f(\xi^{(i)} + s\xi_ie_i) \, ds \right] -
			 			\sum_{i=1}^n \mathbb E \left[ \xi_i^3 \int_0^1 (1-s) \partial_i^{(3)} f(\xi^{(i)} + s\xi_ie_i) \, ds \right].
                \end{align} 
				
			\item If, additionally, $\mathbb E\xi_i^3 =0$ and $\mathbb E\xi_i^4 < \infty$ for all $i\leq n$, 
			then for any $C^4$-smooth function $f:\mathbb R^n \to \mathbb R$ we have
				\begin{align} \label{eq:3-2} \begin{split}
					\mathbb E \left[ Lf(\xi) \right] = \sum_{i=1}^n \mathbb E \left[\xi_i^2 \int_0^1 (1-s) \partial_i^{(4)} f(\xi^{(i)} + s\xi_ie_i) \, ds \right] -
			 			\sum_{i=1}^n \mathbb E \left[\frac{\xi_i^4}{2} \int_0^1 (1-s)^2 \partial_i^{(4)} f(\xi^{(i)} + s\xi_ie_i) \, ds \right],
							\end{split}
				\end{align}
		\end{enumerate}
        where $ \xi^{(i)} := \xi - \xi_i e_i $.
\end{lemma}

\noindent {\it Proof.} We shall demonstrate the proof of \eqref{eq:3-2}. We work similarly for \eqref{eq:3-1}. Recall that $L$ acts as follows:
	\begin{align} \label{eq:1}
		Lf(\xi) = \Delta f(\xi) - \sum_{i=1}^n \xi_i  \partial_i f(\xi).
	\end{align}
Using Taylor's approximation theorem for $\partial_if(\xi)$ \big (in the direction $e_i$ around $0$\big), we get
	\[
	\partial_i f(\xi) - \partial_i f(\xi^{(i)}) 
	= \xi_i \partial^{(2)}_i f(\xi^{(i)}) + \xi^{2}_i \frac{1}{2} \partial^{(3)}_i f(\xi^{(i)}) + 
						\frac{1}{2!} \int_0^{\xi_i} (\xi_i-s)^{2} \partial_i^{(4)} f(\xi^{(i)}+se_i) ds. 	
	\] 
Using the substitution $s \to s\xi_i$, and multiplying both sides by $\xi_i$, we derive
	\[ 
		\xi_i (\partial_if(\xi) - \partial_i f(\xi^{(i)}) ) 
		= \xi_i^{2} \partial_i^{(2)} f(\xi^{(i)}) + \frac{\xi_i^3}{2} \partial_i^{(3)} f(\xi^{(i)})  
				+ \frac{\xi_i^4}{2} \int_0^1 (1-s)^2 \partial_i^{(4)} f(\xi^{(i)} + s \xi_i e_i ) \, ds .
	\]
Now using \eqref{eq:1}, we may write
	\[ 
		L f(\xi)  = \Delta f(\xi) - \sum_{i=1}^n \xi_i (\partial_i f(\xi) - \partial_i f(\xi^{(i)}) ) - 
			\sum_{i=1}^n \xi_i \partial_i f(\xi^{(i)}).
	\]
Taking expectations on both sides, and using the fact that $\xi_i$ are independent \big (in particular $\xi_i$ is independent of $\xi^{(i)}$\big) and satisfy the assumptions of (2), we obtain
	\begin{align} \label{eq-2}
		\mathbb E [ L f(\xi) ] 
			= \mathbb E \left[  \hspace{1mm} \sum_{i=1}^n ( \partial_i^{(2)} f(\xi) - \partial_i^{(2)} f(\xi^{(i)}) ) -
					\frac{1}{2} \sum_{i=1}^n\xi_i^4 \int_0^1 (1-s)^2 \partial_i^{(4)} f (\xi^{(i)} + s \xi_i e_i ) \, ds \right].
	\end{align}
Notice that using the Taylor's approximation theorem for $\partial_i^{(2)} f(\xi)$ \big(and the same substitution as previously used\big), we obtain
	\[ 
		\partial_i^{(2)} f(\xi) - \partial_i^{(2)} f(\xi^{(i)}) 
			= \xi_i \partial_i^{(3)} f(\xi^{(i)}) + \int_0^{\xi_i} (\xi_i - s) \partial_i^{(4)} f (\xi^{(i)} + s e_i) \, ds.
	\]
Plugging this into \eqref{eq-2}, we get
	\begin{align*}
		\mathbb E [L f(\xi) ] &= \sum_{i=1}^n \mathbb E \left [ \xi_i^2\int_0^1 (1-s) \partial_i^{(4)} f(\xi^{(i)} 
					+ s \xi_i e_i) \, ds - \frac{\xi_i^4}{2} \int_0^1 (1-s)^2 \partial_i^{(4)} f(\xi^{(i)} + s \xi_i e_i ) \, ds \right],
	\end{align*} as claimed. The proof is complete. \prend \\

\begin{remark}
    Notice that in the above proof, we used Taylor's approximation theorem twice \big(of order two for $ \partial_i f(\xi)$ and of order one for $ \partial_i^{(2)} f(\xi)$ respectively\big). To prove \eqref{eq:3-1}, one employs Taylor's approximation theorem twice as well \big(of order one for $ \partial_i f(\xi)$ and of order zero for $ \partial_i^{(2)} f(\xi)$ respectively\big). Other details and computations are carried out analogously. 
\end{remark}

\begin{remark}
    We note in passing that one could use Lemma \ref{lem:s-form} to derive universality properties for various problems in communications, statistical learning, and random matrix theory by appropriately choosing the function $f$ in Lemma \ref{lem:s-form} as in \cite{korada2011applications}. \\
\end{remark}

\noindent Using the integral representation in Lemma \ref{lem:s-form} we may already get a (weaker) error approximation bound. Namely, this bound will depend on the following (geometric) parameter\footnote{Note that if $T$ is a compact set with $0\in {\rm int}T$, then $x\mapsto \max_{t\in T}|\langle x, t \rangle |$ is a norm on $\mathbb R^n$-let us denote by $\|\cdot \|$ this norm, and by $X_T$ the corresponding normed space. Thereby, the functional $(x_1, \ldots, x_n) \mapsto \left( \sum_{i=1}^n \|x_i\|^p \right)^{1/p}$ is the induced norm on the $n$th fold $\ell_p$-sum of the space $X_T$, hence the notation $\| \cdot \|_{\ell_p(T)}$.} associated with the set $T$,
	\[
		\|(e_1,\ldots, e_n)\|_{\ell_p(T)}: = \left(\sum_{i=1}^n \max_{t\in T}|\langle t, e_i\rangle|^p \right)^{1/p},
	\] 
for $p \geq 1 $. 

\begin{proposition} \label{prop:weak-1}
Let $T$ be a bounded subset of $\mathbb R^n$. Let $\xi=(\xi_1, \ldots, \xi_n)$ be a random vector with independent coordinates 
such that $\mathbb E \xi_i= 0$ and $\mathbb E \xi_i^2=1$. We have the following cases:
	\begin{enumerate}
		\item If $\sigma_3^3 := \max_{i\leq n} \mathbb E |\xi_i|^3<\infty$, then for all $\beta>0$ we have
			\begin{align}
				|\mathbb EF_\beta(\xi ) - \mathbb EF_\beta(G)| \leq C \beta^2 \sigma_3^3 \|(e_1, \ldots,e_n)\|_{\ell_3(T)}^3.			
			\end{align}
		\item Furthermore, if $\mathbb E\xi_i^3=0$ for all $i\leq n$ and $\sigma_4^4 :=\max_{i\leq n} \mathbb E\xi_i^4 < \infty$, 
		then, for all $\beta>0$ we have
 			\begin{align}
				|\mathbb EF_\beta(\xi ) - \mathbb EF_\beta(G)| \leq C\beta^3 \sigma_4^4 \|(e_1, \ldots,e_n)\|_{\ell_4(T)}^4,
			\end{align}
		where $C>0$ is a universal constant.
	\end{enumerate}
\end{proposition}

\noindent {\it Proof.} Using \eqref{eq:dtc} we may write the following:
	\[
		\mathbb EF_\beta(\xi) - \mathbb EF_\beta(G) = - \mathbb E_\xi \left[ {\mathscr P} LF_\beta (\xi)  \right]
													= - \mathbb E_\xi \left[ L {\mathscr P} F_\beta(\xi) \right].
	\] 
Applying Lemma \ref{lem:s-form} for ``$f= {\mathscr P}F_\beta$'', we infer
	\begin{align} \label{eq:m-1} 
        \begin{split}
		\mathbb EF_\beta(\xi) - \mathbb EF_\beta(G) =	
  -\sum_{i=1}^n \mathbb E \left[\xi_i^2 \int_0^1  (1-s) \partial_i^{(4)}  {\mathscr P}F_\beta (\xi^{(i)} + s\xi_ie_i) \, ds \right] +
			 		  \sum_{i=1}^n  \mathbb E  \left[\frac{\xi_i^4}{2} \int_0^1 (1-s)^2 \partial_i^{(4)} {\mathscr P}F_\beta (\xi^{(i)} + s\xi_ie_i) \, ds \right] .
        \end{split}
	\end{align} 
In view of \eqref{eq:partial-pt} we have
	\begin{align} \label{eq:m-2}
		\partial_i^{(4)}{\mathscr P} F_\beta(\xi^{(i)} + s\xi_i e_i)  = \int_0^\infty e^{-4t} P_t (\partial_i^{(4)}F_\beta)(\xi^{(i)} + s\xi_i e_i) \, dt.
	\end{align} 
Note that 
	\begin{align} \label{eq:m-3}
		P_t (\partial_i^{(4)}F_\beta)(\xi^{(i)} + s\xi_i e_i)  = 
			\mathbb E_G \left[ \partial_i^{(4)}F_\beta \left(e^{-t}(\xi^{(i)} + s\xi_i e_i) +  \sqrt{1-e^{-2t}} G\right) \right].
	\end{align} 
From Lemma \ref{lem:bd-der-Fbeta} we conclude that
	\begin{align}
		|\partial_i^{(4)}F_\beta \left(e^{-t}(\xi^{(i)} + s\xi_i e_i) +  \sqrt{1-e^{-2t}} G\right)| \leq 26\beta^3 \max_{t\in T}|t_i|^4,
	\end{align} 
a.s. Therefore, combining the latter pointwise estimate with all the above we infer that
	\begin{align}
		\left| \mathbb EF_\beta(\xi) - \mathbb EF_\beta(G) \right| \leq \frac{26\beta^3}{4}  \sum_{i=1}^n \mathbb E \xi_i^2 \max_{t\in T}|t_i|^4 + 
			\frac{26\beta^3}{8} \sum_{i=1}^n \mathbb E \xi_i^4 \max_{t\in T} |t_i|^4 \leq \frac{39\beta^3}{4} \sigma_4^4 \sum_{i=1}^n \max_{t\in T}|t_i|^4.
	\end{align}
The proof is complete. \prend \\

\begin{corollary} \label{corr-main}
	Let $T$ be a finite subset of $\mathbb R^n$ and let $\xi= (\xi_1, \ldots, \xi_n)$ be a random vector with independent coordinates 
	such that $\mathbb E \xi_i=0$ and $\mathbb E\xi^2 =1$ for all $i\leq n$. We have the following assertions:
		\begin{enumerate}
			\item If $\sigma_3^3 : = \max_{i\leq n} \mathbb E |\xi_i|^3 < \infty$, then  
				\begin{align}
					\left| \mathbb E\sup_{t\in T} \langle \xi, t\rangle - g(T) \right| \leq C \sigma_3 \|(e_i)_{i\leq n}\|_{\ell_3(T)} (\log |T|)^{2/3}. 				
				\end{align}
			\item Moreover, if $\mathbb E\xi_i^3=0$ and $\sigma_4^4 := \max_{i\leq n} \mathbb E |\xi_i|^4 <\infty$, then 
				\begin{align}
					\left| \mathbb E\sup_{t\in T} \langle \xi, t\rangle - g(T) \right| \leq C \sigma_4 \|(e_i)_{i\leq n}\|_{\ell_4(T)} (\log |T|)^{3/4},
				\end{align}
		\end{enumerate}
		where $C>0$ is a universal constant.
\end{corollary}

\noindent {\it Proof.} We only prove (2). By virtue of Lemma \ref{lem:prop-Fbeta} we have
	\[
		\left| \mathbb E\sup_{t\in T} \langle \xi, t\rangle - g(T) \right| \leq \frac{\log|T|}{\beta} + \left| \mathbb E F_\beta (\xi) - \mathbb E F_\beta(G)\right|.
	\] 
	Invoking Proposition \ref{prop:weak-1} we obtain
	\[
		\left| \mathbb E\sup_{t\in T} \langle \xi, t\rangle - g(T) \right|  \leq \frac{\log |T|}{\beta} + C\beta^3 \sigma_4^4 \|(e_i)_{i\leq n}\|_{\ell_4(T)}^4,
	\] 
	for all $\beta>0$. We optimize with respect to $\beta$ by choosing $\beta= \frac{(\log |T|)^{1/4}}{\sigma_4 \|(e_i)_{i\leq n}\|_{\ell_4(T)}}$ . \prend \\

\noindent For establishing the announced (stronger) error approximation bound in terms of $R_4(T)$, 
we assume additionally that the random variables are uniformly bounded. Under
this assumption one may bound in a more efficient way the $\mathbb E_{\mu_x}[\ell_i^4]$, and thereby the $|\partial^{(4)} F_\beta(x)|$. Toward this end, we shall need the following technical result which also appeared in \cite{Tal-Gibbs}. Since its proof is short we reproduce it here for reader's convenience.

\begin{lemma} \label{lem:lip-log-avg}
	Let $T\subset \mathbb R^n$ be bounded and let $\beta>0$. For fixed $ 1\leq i \leq n$, any $x=(x_1,\ldots, x_i , \ldots x_n)\in \mathbb R^n$ and 
	$x' =(x_1, \ldots, x_i', \ldots, x_n)$ we have
	\begin{align}
		\mathbb E_{\mu_x}[ \ell_i^4] \leq e^{2\beta R_\infty |x_i-x_i'|} \mathbb E_{\mu_{x'} }[ \ell_i^4],
	\end{align}
         where $R_\infty := \sup_{t \in T} \norm{t}_{\infty}$ and $ \ell_i := \ell_i(t)$.
\end{lemma}

\noindent {\it Proof.} Let $\psi (z)  = \mathbb E_{\mu_z} [\ell_i^4]$. Applying the Fundamental Theorem of Calculus we obtain
\begin{align} \label{eq:ftclog}
\log \psi(x) - \log\psi(x') = \int_{x_i'}^{x_i} (\partial_i \log \psi) (x^{(i)} + se_i) \, ds, 
\end{align} 
where $x^{(i)}= (x_1, \ldots, x_{i-1},0,x_{i+1}, \ldots, x_n)$. On the other hand
\[
\partial_i \psi(z) = \beta \mathbb E_{\mu_{z}}[\ell_i^5] - \beta \mathbb E_{\mu_{z}}[ \ell_i^4] \cdot \mathbb E_{\mu_{z}} [\ell_i ],
\] 
for all $z$. Thereby, we get
\[
|\partial_i \psi(z)| \leq 2\beta R_\infty \mathbb E_{\mu_{z}} [  \ell_i^4 ] = 2\beta R_\infty \psi(z),
\]
which implies that 
\[
\abs{\partial_i \log( \psi(z))} \leq 2 \beta R_\infty,
\]
for all $z \in \mathbb R^n $.
Combing the latter estimate with the identity \eqref{eq:ftclog} the result follows. \prend \\

\noindent The next result is the counterpart of Proposition \ref{prop:weak-1} in the case of uniformly bounded random variables. More precisely, we have the 
following.

\begin{proposition} \label{prop:strong-1}
	Let $T\subset \mathbb R^n$ be bounded and let $\beta>0$. Then, for any random vector $\xi=(\xi_1, \ldots, \xi_n)$ with independent coordinates, 
	where $\mathbb E\xi_i = \mathbb E\xi_i^3=0$, $\mathbb E\xi_i^2=1$ and $|\xi_i| \leq M$ a.s. for all $i\leq n$, we have that
	\begin{align}
		\left|\mathbb E F_\beta(\xi) - \mathbb EF_\beta(G) \right| \leq C \beta^3 M^4 R_4^4 e^{2 \beta M R_\infty}.
	\end{align}
\end{proposition}

\noindent {\it Proof.} As in the proof of Proposition \ref{prop:weak-1}, using \eqref{eq:dtc} we may write the following:
	\[
		\mathbb EF_\beta(\xi) - \mathbb EF_\beta(G) = - \mathbb E_\xi \left[ {\mathscr P} LF_\beta (\xi)  \right]
													= - \mathbb E_\xi \left[ L {\mathscr P} F_\beta(\xi) \right].
	\] 
Applying Lemma \ref{lem:s-form} for ``$f= {\mathscr P}F_\beta$'', we infer
	\begin{align} \label{eq:m-1} \begin{split}
		\mathbb EF_\beta(\xi) - \mathbb EF_\beta(G) =	-			
				 \sum_{i=1}^n \mathbb E \left[\xi_i^2 \int_0^1  (1-s) \partial_i^{(4)}  {\mathscr P}F_\beta (\xi^{(i)} + s\xi_ie_i) \, ds \right] +
			 		  \sum_{i=1}^n  \mathbb E  \left[\frac{\xi_i^4}{2} \int_0^1 (1-s)^2 \partial_i^{(4)} {\mathscr P}F_\beta (\xi^{(i)} + s\xi_ie_i) \, ds \right] .
			\end{split}
	\end{align} 
In view of \eqref{eq:partial-pt} we have
	\begin{align} \label{eq:m-2}
		\partial_i^{(4)}{\mathscr P} F_\beta(\xi^{(i)} + s\xi_i e_i)  = \int_0^\infty e^{-4t} P_t (\partial_i^{(4)}F_\beta)(\xi^{(i)} + s\xi_i e_i) \, dt.
	\end{align} 
Note that 
	\begin{align} \label{eq:m-3}
		P_t (\partial_i^{(4)}F_\beta)(\xi^{(i)} + s\xi_i e_i)  = 
			\mathbb E_G \left[ \partial_i^{(4)}F_\beta \left(e^{-t}(\xi^{(i)} + s\xi_i e_i) +  \sqrt{1-e^{-2t}} G\right) \right].
	\end{align} 
Recall that Lemma \ref{lem:bd-der-Fbeta} yields the (pointwise) bound 
	\begin{align}	\label{eq:m-4a}
		\left| \partial_i^{(4)}F_\beta \left(e^{-t}(\xi^{(i)} + s\xi_i e_i) +  \sqrt{1-e^{-2t}} G\right) \right|  \leq 26 \beta^3 \mathbb E_{\mu_{z_i}}[\ell_i^4],
	\end{align} 
a.s., where $z_i= e^{-t}(\xi^{(i)} + s\xi_i e_i) +  \sqrt{1-e^{-2t}} G$. \\
If we set $z= e^{-t}\xi + \sqrt{1-e^{-2t} }G$, an appeal to 
Lemma \ref{lem:lip-log-avg} (applied for ``$x' =z$'' and ``$x=z_i$'') yields the (almost sure) bound
	\begin{align} \label{eq:m-4b}
		\mathbb E_{\mu_{z_i}}[\ell_i^4] \leq e^{2\beta (1-s) R_\infty e^{-t}|\xi_i|} \mathbb E_{\mu_{z} } [\ell_i^4] 
														\leq e^{2\beta  M R_\infty } \mathbb E_{\mu_{z}} [\ell_i^4].
	\end{align}
Thus, \eqref{eq:m-4a} and \eqref{eq:m-4b} imply that
	\begin{align} \label{eq:m-5a}
		\left| \partial_i^{(4)}F_\beta \left(e^{-t}(\xi^{(i)} + s\xi_i e_i) +  \sqrt{1-e^{-2t}} G\right) \right|  
				\leq 26 \beta^3 e^{2\beta M R_\infty} \mathbb E_{\mu_z}[\ell_i^4],
	\end{align} 
for all $\xi$, and for all $0\leq s\leq 1$ and $t\geq 0$. Taking into account \eqref{eq:m-2}, \eqref{eq:m-3}, and \eqref{eq:m-5a} we obtain
	\begin{align} \label{eq:m-5}
		|\partial_i^{(4)}{\mathscr P}F_\beta (\xi^{(i)} + s \xi_i e_i)| \leq 
				26\beta^3 e^{2\beta M R_\infty} \int_0^\infty e^{-4t} \mathbb E_G \mathbb E_{\mu_{e^{-t} \xi + \sqrt{1-e^{-2t}}G } }[\ell_i^4] \, dt.
	\end{align}
Finally, combining \eqref{eq:m-5} with \eqref{eq:m-1} we obtain
	\begin{align*}
		\left| \mathbb EF_\beta(\xi) - \mathbb EF_\beta(G) \right| &\leq 26 \beta^3 M^3e^{2\beta M R_\infty} 				
				  \int_0^1 (1-s) \left( \int_0^\infty e^{-4t} \mathbb E_{\xi,G} \mathbb E_{\mu_{e^{-t} \xi + \sqrt{1-e^{-2t}} G} } \left[\Sigma_{i = 1}^n \ell_i^4 \right] \, dt\right) \, ds \\
& \hfill +\frac{26\beta^3}{2}M^4 e^{2\beta M R_\infty} \int_0^1 (1-s)^2 \left( \int_0^\infty e^{-4t} \mathbb E_{\xi,G} \mathbb E_{\mu_{e^{-t} \xi + \sqrt{1-e^{-2t}} G} } \left[\Sigma_{i = 1}^n \ell_i^4 \right] \, dt\right) \, ds \\
				  & \leq 5 \beta^3  M^4 R_4^4 e^{2\beta M R_\infty},
	\end{align*}
where we have used the fact that $M\geq (\mathbb E \xi_i^2)^{1/2}=1$. The proof is complete. \prend \\

\begin{remark} In estimate \eqref{eq:m-4b}, instead of using Lemma \ref{lem:lip-log-avg}, one might hope to show that there is a ``majorant" probability measure $\nu$ (independent of $i$) on $T$ such that
\[
\mathbb E_{ \mu_{z_i}} [\ell_i^4] \leq C \, \mathbb E_{ \nu } [\ell_i^4]
\]
for all $i \in \{1,2,\ldots,n \}$ and for some universal constant $C > 0$, thereby implying that
\begin{align} \label{for contra}
\sum_{i=1}^n \mathbb E_{ \mu_{z_i}} [\ell_i^4] \leq \sum_{i=1}^n C \, \mathbb E_{ \nu } [\ell_i^4] \leq C \, R_4^4,
\end{align}
and hence sharpening Proposition \ref{prop:strong-1} (by removing exponential dependence on $R_\infty$). We next show that such a ``majorant" measure $\nu$ does not exist in general.

\medskip

\noindent First notice that in estimate \eqref{eq:m-4b}, one may assume without loss of generality that $ z_i = \xi^{(i)} + s\xi_i e_i $. Also, we see that in \eqref{for contra}, right-hand side does not depend on $\xi$ nor on $\beta > 0$. Recall that $ s \in [0,1]$. For $\theta>0$, consider the set $T_\theta = \big \{ - \theta e_j \, | \, 1 \leq j \leq n \}$, where $\{e_j\}_{j=1}^n$ is the canonical basis in $\mathbb R^n $  and assume $\xi_i$ are $i.i.d.$ Rademacher random variables . We get
\begin{align*}
 \mathbb E_{ \mu_{\xi^{(i)} + s\xi_i e_i}} [\ell_i^4] &= \sum_{j=1}^n [\ell_i (- \theta e_j)]^4 \, \, \, \frac{e^{\beta \langle -\theta e_j , \,  \xi^{(i)} + s\xi_ie_i \rangle}}{\sum_{j=1}^ne^{\beta \langle -\theta e_j , \,  \xi^{(i)} + s\xi_ie_i \rangle}  } \\
 &= \theta^4 \, \, \, \frac{e^{- \beta \theta s \xi_i }}{e^{- \beta \theta \xi_1} + .... + e^{- \beta \theta s \xi_i} + .... + e^{- \beta \theta \xi_n}} \\
 &= \theta^4 \frac{e^{-\theta/2}}{e^{-\theta}+....+e^{-\theta/2}+...+e^{-\theta}}
\end{align*}
for $\xi_i =1$, $\beta=1$, and $s=1/2$.

\medskip

\noindent Plugging the above in \eqref{for contra} and choosing $\theta$ large enough, we get
\begin{align*}
\sum_{i=1}^n \mathbb E_{ \mu_{\xi^{(i)} + s\xi_i e_i}} [\ell_i^4] \asymp n \theta^4 = n R_4^4(T_\theta), 
\end{align*}
\noindent and hence \eqref{for contra} cannot hold for any universal constant $C>0$. \\
\end{remark}

\noindent We are now ready to prove our main result: 

\begin{theorem} \label{thm-main}
Let $T\subset \mathbb R^n$ be a finite set. Then, for any random vector $\xi=(\xi_1,\ldots, \xi_n)$ with independent coordinates where $\mathbb E\xi_i = \mathbb E\xi_i^3=0$, $\mathbb E\xi_i^2=1$, and $| \xi | \leq M$ a.s. we have
	\begin{align}
		\left| \mathbb E \sup_{t\in T} \langle \xi, t\rangle - g(T) \right| \leq C \max \left\{ MR_4(\log |T|)^{3/4}, MR_\infty \log |T| \right\},
	\end{align} where $C>0$ is a universal constant.
\end{theorem}

\noindent {\it Proof.} By virtue of Lemma \ref{lem:prop-Fbeta} and Proposition \ref{prop:strong-1} we have
	\[
		\left| \mathbb E \sup_{t\in T} \langle \xi, t\rangle - g(T) \right| \leq \frac{\log|T|}{\beta} + 4 \beta^3 M^4 R_4^4 e^{2 \beta MR_\infty},
	\] 
for any $\beta>0$. We choose 
	\[
	\beta = \min \left \{ (M R_\infty)^{-1}, \frac{(\log |T|)^{1/4}}{M R_4} \right \}
	\] 
to conclude. \prend \\

\begin{remark}
One may modify Lemma \ref{lem:lip-log-avg} to obtain a similar bound for $ \mathbb E_{\mu_x} [\ell_i^3] $. Then, by combining this bound with estimate \eqref{eq:3-1}, we obtain a result analogous to Preposition \ref{prop:strong-1} in terms of $R_3$. Finally, one may dispose the assumption $\mathbb E \xi_{i}^3 = 0$ in Theorem \ref{thm-main} and obtain the following weaker comparison estimate:
    \[
    \left| \mathbb E \sup_{t\in T} \langle \xi, t\rangle - g(T) \right| \leq C \max \left\{ MR_3(\log |T|)^{2/3}, MR_\infty \log |T| \right\},
    \]
    where $C>0$ is a universal constant.   
\end{remark}


\medskip

\section{Discussions and further remarks}

\subsection{Discussion: Theorem \ref{thm:1-1} versus Theorem \ref{thm-main} } \label{sec-4:picture}

\begin{figure}[ht!] 
\includegraphics[width=5in]{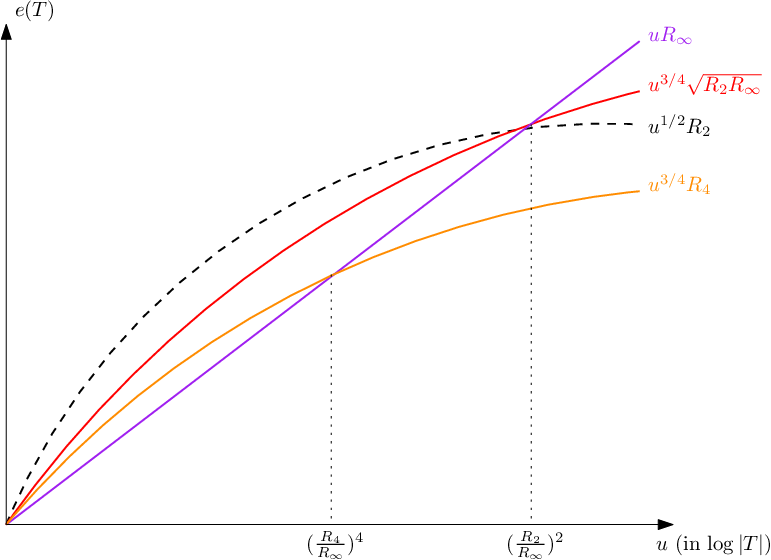}
\centering
\caption[error window]{Phase transition in error bound}
\label{error-fig}
\end{figure}

\noindent In this subsection we discuss the refinement on the proximity bound 
\[
	e(T): = |r(T) - g(T)|, 
\]
where $r(T)$ and $g(T)$ are Rademacher and Gaussian complexities of T, respectively.\\

\noindent Let us define the following curves (associated with $T$):
\begin{align}
k(u) = u^{1/2}R_2, \quad {\color{red} t(u) = u^{3/4}\sqrt{R_2R_\infty} }, \quad {\color{orange} s_1(u)=u^{3/4}R_4} , \quad {\color{purple} s_2(u) = uR_\infty},
\end{align}
where these are respectively the trivial bound, Talagrand's bound and the bounds from Theorem \ref{thm-main}.

\noindent We shall introduce the following threshold values:
\begin{align}
u_1\equiv u_1(T):= (R_4/R_\infty)^4, \quad u_2 \equiv u_2(T) :=(R_2/R_\infty)^2.
\end{align}

\noindent Note that for any bounded set $T\subset \mathbb R^n$ we have $u_1(T) \leq u_2(T) \leq n$. \\

\noindent Recall that due to the fact that both processes $ \{\langle \mathcal{E} ,t \rangle \}_{t\in T}$ and $ \{ \langle G, t\rangle \}_{t\in T}$ are sub-gaussian  
with constants approximately equal to $R_2(T)$ \big(see Remark \ref{rem:1-2}\big) we readily get the (trivial) bound. 
\begin{align*}
|r(T)-g(T)| \leq CR_2\sqrt{\log |T|}.
\end{align*}

\noindent Therefore, the error bound $t(u)$, ($u=\log |T|$) offered by Theorem \ref{thm:1-1} improves upon the trivial one when $u \ll u_2$. In view of $R_4\leq \sqrt{R_2 R_\infty}$, 
one may check that in this range we have $s(u):=\max\{s_1(u),s_2(u)\} \leq t(u)$. 
More precisely, if $R_4 =\delta \sqrt{R_2 R_\infty}$, with $\delta \ll 1$ the window $[u_1,u_2]= [\delta^4u_2, u_2]$ of phase transition starts appearing where
	\[
		s(u) =\begin{cases} s_1(u), & u\leq u_1 \\
						 s_2(u), & u_1< u \leq u_2.
						\end{cases}
	\]

\noindent This refined error bound is illustrated in the Figure \ref{error-fig} above. Note that this is a dynamic phenomenon which depends on the index set $T$ and fades down as 
$u_1(T) \to u_2(T)$. However, the relative position of the aforementioned curves/parameters persists for each $T$. \\

\begin{remark} As in Theorem \ref{thm:1-1} we always have (due to sub-gaussian behavior \cite[Chapter 7]{vershynin2018high}) that $\abs{g(T)}$, $\abs{\mathbb E \sup_{t \in T} \langle \xi, t \rangle}$ $\leq$ $C M R_2 \sqrt{\log |T|}$ , and hence Theorem \ref{thm-main} is of interest only when its right-hand side is smaller than this quantity, that is, when
    \[
     R_4 \big (\log |T| \big )^{1/4} \, , \,  R_\infty \big (\log |T| \big )^{1/2} \leq R_2.
    \]
\end{remark}

\begin{remark}
Notice that one cannot hope to generalize Talagrand's Theorem \ref{thm:1-1} to include $\xi_i$ with sub-exponential tails. To see this, choose $T= \big \{ e_i \mid 1\leq i \leq n \big \} $, that is, the canonical basis in $\mathbb R^n$, and $\xi_i$ to be $i.i.d.$ Laplace (two-sided exponential) distribution, that is, each $\xi_i$ has the density $\frac{1}{2} e^{- \abs{x}}$ where $x \in \mathbb R $ \,. Then, the left-hand side of Theorem \ref{thm:1-1} is of the order $\log n $ whereas the right-hand side is of the order $ {(\log n)}^{3/4} $.
\end{remark}

\subsubsection{Examples} \label{sec-4:examples}

In this subsection we construct a family of finite sets $T \subset \mathbb R^n$ for a wide spectrum of cardinalities with the additional property $(R_4 / R_\infty)^4 \ll ( R_2 / R_\infty)^2$. \\

\noindent Consider any $A \subset \{ -1, 1 \}^n$ such that $ |A| = 2^k$, where $k \ll n$. For a
decreasing sequence of positive numbers $d_1> d_2 > \ldots > d_n>0$ we consider the diagonal matrix 
	\[
		D : = {\rm diag} (d_1, \ldots, d_n) = \begin{pmatrix} d_1 & & & \\
												& d_2 & &\\ 
													& & \ddots &\\
													& & & d_n
									\end{pmatrix}. 
	\] \\

\noindent Write $d$ for the vector $(d_1, \ldots, d_n)$ and set $T_d := D(A) = \{ Da \,\,|\,\, a \in A \}$, and notice that, 
for every $t \in T_d$, $t_i = \pm d_i $, and so $\abs{t_i}=d_i$. Therefore, we get $\|t\|_p= \|d\|_p$ for all $p>0$, and thus
	\begin{align}
		u_1(T_d) = (R_4/R_\infty)^4 = (\|d\|_4 / \|d\|_\infty)^4, \\
		u_2(T_d) = (R_2/R_\infty)^2 = (\|d\|_2/ \|d\|_\infty)^2,
	\end{align}
and $\log |T_d| \asymp k$. Taking into account the above features, we may construct subsets $T_d \subset \mathbb R^n $ with the property that
$\log |T_d| < (R_2 / R_\infty)^2$ and $(R_4/ R_\infty)^4 \ll (R_2/ R_\infty)^2$. For example, fix $k< \sqrt{n}$ and take $d_j \sim j^{-1/4}$ for $ j\leq n$.

\subsection{Further remarks}

\subsubsection{Universality} \label{sec-4:univ}

Talagrand in \cite{Tal-Gibbs} proved that the energy of the ground state of the Sherrington-Kirkpatrick's (SK) model does not change when the Gaussian interactions $g_{ij}$ are replaced by Rademacher (symmetric Bernoulli) $\eta_{ij}$ interactions. In particular, he proved the following: 

\begin{corollary} [Talagrand]
    \label{cor:4-1} Consider i.i.d. standard Gaussian random variables $(g_{ij})_{1\leq i<j \leq N}$ and i.i.d. Rademacher random variables $(\eta_{ij})_{1\leq i<j \leq N}$. Then 
    \[
    \abs{\frac{1}{N^{3/2}} \mathbb E \, \sup \sum_{1 \leq i < j \leq N }  g_{ij} \sigma_i \sigma_j - \frac{1}{N^{3/2}} \mathbb E \, \sup \sum_{1 \leq i < j \leq N} \eta_{ij} \sigma_i \sigma_j } \leq C N^{-1/4},
    \]
    where the supremum is over $ \sigma = (\sigma_1, \sigma_2, ...., \sigma_N) \in \{-1, 1 \}^N $. \\
\end{corollary}

\noindent Here, each $\sigma$ is interpreted as a configuration of an $N$-spin physical system. Above Corollary \ref{cor:4-1} was later generalized by Carmona and Hu in \cite{CH-univ} to other interactions (environments) beyond Rademacher interactions $\eta_{ij}$.

\begin{corollary} [Carmona and Hu]
    \label{corr-4-2} Let $(g_{ij})_{1 \leq i<j \leq N}$ be i.i.d. standard Gaussian random variables, and let $(\eta_{ij})_{1 \leq i<j \leq N}$ be independent with $\mathbb E \eta_{ij} =0$ and $\mathbb E \eta_{ij}^2=1$. We have the following assertions:
    \begin{enumerate}
        \item If $ \mu_3^3 := \max_{1 \leq i<j \leq N} \mathbb E \abs{\eta_{ij}}^3 < \infty$, then 
        \begin{align}
             \abs{\frac{1}{N^{3/2}} \mathbb E \, \sup \sum_{1 \leq i<j \leq N }  g_{ij} \sigma_i \sigma_j - \frac{1}{N^{3/2}} \mathbb E \, \sup \sum_{1 \leq i<j \leq N} \eta_{ij} \sigma_i \sigma_j } \leq C \mu_3 N^{-1/6}
        \end{align}
        \item Moreover, if $ \mathbb E \eta_{ij}^3=0 $ and $ \mu_4^4 := \max_{1 \leq i<j \leq N} \mathbb E \abs{\eta_{ij}^4} < \infty $, then 
        \begin{align}
          \abs{\frac{1}{N^{3/2}} \mathbb E \, \sup \sum_{1 \leq i<j \leq N }  g_{ij} \sigma_i \sigma_j - \frac{1}{N^{3/2}} \mathbb E \, \sup \sum_{1 \leq i<j \leq N} \eta_{ij} \sigma_i \sigma_j } \leq C \mu_4 N^{-1/4}   
        \end{align}
    \end{enumerate}
    where $C>0$ is a universal constant. \\
\end{corollary}

\noindent We next show that the above results are special cases of a more general universality principle for the following  \textit{random order-m} tensors:

\begin{corollary} \label{corr-tensor}
    Let $\sigma = (\sigma_1, \sigma_2, \dots, \sigma_N) \in \{-1,1\}^N \subset \mathbb R^N $ and define a random order-m tensor as follows:
    \[
    (i_1, i_2, \dots, i_m)  \mapsto \eta_{i_1 i_2 ....i_m} \prod_{j=1}^m\sigma_{i_j} 
    \]
    where $ 1 \leq i_1 < i_2 < \dots < i_m \leq N $, $ 1 < m < (N-1)$, and $\eta_{i_1 i_2 ....i_m}$ are independent with $\mathbb E \, \eta_{i_1 i_2 ....i_m} = 0 $ and $\mathbb E \, \eta^{2}_{i_1 i_2 ....i_m} = 1 $. Let $g_{i_1 i_2 ....i_m}$ be $i.i.d.$ standard Gaussian random variables.
    \begin{enumerate}
        \item If $ \mu_3^3 := \max_{1 \leq i_1 < i_2 < \dots < i_m \leq N} \mathbb E \abs{\eta_{i_1 i_2 ....i_m}}^3 < \infty $, then
        \begin{align*}
        &\abs{ \frac{1}{{\binom{N}{m}}^{1/2} N^{1/2}} \mathbb E \sup_\sigma \sum_{1 \leq i_1 < i_2 < \dots < i_m \leq N} \eta_{i_1 i_2 ....i_m} \prod_{j=1}^m\sigma_{i_j} - \frac{1}{{\binom{N}{m}}^{1/2} N^{1/2}}  \mathbb E \sup_\sigma \sum_{1 \leq i_1 < i_2 < \dots < i_m \leq N} g_{i_1 i_2 ....i_m} \prod_{j=1}^m\sigma_{i_j}} \\
        &\leq C \mu_3 \Bigg ( {\frac{N}{\binom{N}{m}}} \Bigg )^{1/6}.
        \end{align*}
        \item Moreover, if $ \mathbb E \eta_{i_1 i_2 ....i_m}^3 = 0 $, and $\mu_4^4 := \max_{1 \leq i_1 < i_2 < \dots < i_m \leq N} \mathbb E \abs{\eta_{i_1 i_2 ....i_m}}^4 < \infty $, then 
        \begin{align*}
        &\abs{ \frac{1}{{\binom{N}{m}}^{1/2} N^{1/2}} \mathbb E \sup_\sigma \sum_{1 \leq i_1 < i_2 < \dots < i_m \leq N} \eta_{i_1 i_2 ....i_m} \prod_{j=1}^m\sigma_{i_j} - \frac{1}{{\binom{N}{m}}^{1/2} N^{1/2}}  \mathbb E \sup_\sigma \sum_{1 \leq i_1 < i_2 < \dots < i_m \leq N} g_{i_1 i_2 ....i_m} \prod_{j=1}^m\sigma_{i_j}} \\
        &\leq C \mu_4 \Bigg ( {\frac{N}{\binom{N}{m}}} \Bigg )^{1/4},
        \end{align*}
    \end{enumerate}
    where $C>0$ is a universal constant. 
\end{corollary}

\noindent {\it Proof.} Use Corollary \ref{corr-main} with $T= \Bigg \{ \bigg(\frac{\prod_{j=1}^m\sigma_{i_j} }{{\binom{N}{m}}^{1/2} N^{1/2}}\bigg)_{1 \leq i_1 < i_2 < \dots < i_m \leq N} \,\,\, \Big| \,\,\,  \sigma = (\sigma_1,....,\sigma_N) \in \{-1,1\}^N  \Bigg \} $. \\

\begin{remark}
    The normalization used in the above Corollary \ref{corr-tensor} is proper one and recovers Talagrand's normalization used in \cite{Tal-Gibbs} for $m = 2$. This can be seen from the fact that the expected supremum of the Gaussian tensor normalized in this way is bounded away from zero and infinity. More precisely,
    \begin{align*}
       0 < c_m <  \frac{1}{{\binom{N}{m}}^{1/2} N^{1/2}}  \mathbb E \sup_\sigma \sum_{1 \leq i_1 < i_2 < \dots < i_m \leq N} g_{i_1 i_2 ....i_m} \prod_{j=1}^m\sigma_{i_j} < C < \infty, 
    \end{align*}
where constant $c_m > 0 $ depends on $m$, but it is independent of the dimension $N$, and $C$ is as usual, a universal constant. The above can be proved using a basic combinatorial argument integrated with Sudakov Minoration for Gaussian processes. We leave the details to the interested reader. We refer the reader to \cite[Chapter 2]{Tal-spin} for a physical interpretation of the above result when $m = 2$ in the context of spin glasses.  
\end{remark}

\subsubsection{A more general setting} \label{sec-4:general}

Let us note that the results in Section \ref{sec-3} can be formulated for an arbitrary measure $\mu$ supported on the index set $T$. To see this note that for a finite set $T$ we may write
	\begin{align*}
		F_\beta(x) = \frac{1}{\beta} \log \left( \sum_{t\in T} e^{\beta \langle x, t\rangle} \right) = 
			\frac{\log |T|}{\beta} + \frac{1}{\beta} \log \left( \int_T e^{\beta {\langle x,t \rangle}} \, d\nu(t) \right),	
	\end{align*}
where $\nu$ is the uniform probability measure on $T$. The latter integral is an appropriate re-scaling of the {\it logarithmic Laplace transform} of $\nu$. Recall that for a probability measure $\mu$ on $\mathbb R^n$ the logarithmic Laplace transform of $\mu$ is defined by
	\begin{align}
		\Lambda_\mu (\xi)  \stackrel{\rm def} = \log \left ( \int_{\mathbb R^n} e^{\langle \xi, x \rangle} \, d \mu(x) \right), \quad \xi\in\mathbb R^n.
	\end{align}

\noindent We refer the reader to \cite{LW} and the references therein for background information on the log-Laplace transform. With this notation the corresponding Gibbs' measure is given by
\begin{align}
    \frac{d \mu_\xi}{d\mu}(t) = \frac{e^{\langle \xi,t\rangle}}{\int e^{\langle \xi,t \rangle} \, d\mu(t)}.
\end{align} 
Thus, differentiation under the integral sign yields
\begin{align}
    \nabla \Lambda_\mu(\xi) = (\mathbb E_{\mu_\xi}[\ell_i])_{i=1}^n, 
    \quad {\rm Hess} \, \Lambda_\mu(\xi) = {\rm Cov}(\mu_\xi).
\end{align} 
See e.g., \cite{K, KM} for details. Furthermore, Lemma \ref{lem:bd-der-Fbeta} becomes as follows:

\begin{lemma} \label{lem:bd-der-general}
    Let $\mu$ be a probability measure supported on the bounded $T\subset \mathbb R^n$. Then, 
    \begin{align}
        \left| \frac{\partial^3}{\partial \xi_i^3}\Lambda_\mu(\xi)\right| \leq 6 \mathbb E_{\mu_\xi}\left [|\ell_i|^3 \right ] \leq 6\max_{t\in T}|t_i|^3, \quad \left| \frac{\partial^4}{\partial \xi_i^4}\Lambda_\mu(\xi)\right| \leq 26 \mathbb E_{\mu_\xi}\left [|\ell_i|^4 \right ] \leq 26\max_{t\in T}|t_i|^4.    
    \end{align}
\end{lemma}

\noindent Using the above Lemma one can prove the following extension of Proposition \ref{prop:weak-1}.
\begin{proposition}
Let $T\subset \mathbb R^n$ be a bounded set and let be $\mu$ be a probability measure on $T$. Let $\xi=(\xi_1, \ldots, \xi_n)$ be a random vector with independent coordinates 
such that $\mathbb E \xi_i= 0$ and $\mathbb E \xi_i^2=1$. We have the following cases:
	\begin{enumerate}
		\item If $\sigma_3^3 := \max_{i\leq n} \mathbb E |\xi_i|^3<\infty$, then we have
			\begin{align}
				\left|\mathbb E \Lambda_\mu( \xi) - \mathbb E \Lambda_\mu( G)\right| \leq C  \sigma_3^3 \|(e_1, \ldots,e_n)\|_{\ell_3(T)}^3.			
			\end{align}
		\item Furthermore, if $\mathbb E\xi_i^3=0$ for all $i\leq n$ and $\sigma_4^4 :=\max_{i\leq n} \mathbb E\xi_i^4 < \infty$, 
		then we have
 			\begin{align}
				\left|\mathbb E \Lambda_\mu( \xi) - \mathbb E \Lambda_\mu( G)\right| \leq C \sigma_4^4 \|(e_1, \ldots,e_n)\|_{\ell_4(T)}^4,
			\end{align}
		where $C>0$ is a universal constant. \\
	\end{enumerate}
\end{proposition}

\noindent One may generalize Lemma \ref{lem:lip-log-avg} to the following: 
\begin{lemma}
        \noindent Let $T \subset \mathbb R^n$ be bounded and let $\mu$ be a probability measure supported on $T$. For any $f:T \to [0,\infty)$ we may define $ \psi(x) = \log \mathbb E_{\mu_x}[f]$. Then, for all $x,y\in \mathbb R^n $, we have
        \[
        |\psi (x) - \psi(y) |\leq \sup_{t\in T}|\langle t, x-y \rangle|. 
        \]
        In particular, $|\psi (x) - \psi(y) | \leq R_\infty (T) \|x-y\|_1$ for all $x,y \in \mathbb R^n$.   \\ 
\end{lemma}

\noindent Using the above lemma along with Lemma \ref{lem:bd-der-general}, one may obtain the following extension of Proposition \ref{prop:strong-1} :
\begin{proposition}
	Let $T\subset \mathbb R^n$ be bounded and let $\mu$ be a probability measure supported on $T$. Then, for any random vector $\xi=(\xi_1, \ldots, \xi_n)$ with independent coordinates, 
	where $\mathbb E\xi_i = \mathbb E\xi_i^3=0$, $\mathbb E\xi_i^2=1$ and $|\xi_i| \leq M$ a.s. for all $i\leq n$, we have that
	\begin{align}
		\left|\mathbb E \Lambda_\mu( \xi) - \mathbb E \Lambda_\mu( G)\right| \leq C M^4 R_4^4(T) e^{ M R_\infty(T)},
	\end{align}
where $C>0$ is a universal constant. 
\end{proposition}


\medskip


\bibliographystyle{alpha}
\bibliography{Suprema}

\begin{thebibliography}{AAGM15}

\bibitem[AAGM15]{AGM-book}
Shiri Artstein-Avidan, Apostolos Giannopoulos, and Vitali~D. Milman.
\newblock {\em Asymptotic geometric analysis. {P}art {I}}, volume 202 of {\em Mathematical Surveys and Monographs}.
\newblock American Mathematical Society, Providence, RI, 2015.

\bibitem[ALMT14]{ALMT_edge}
Dennis Amelunxen, Martin Lotz, Michael~B McCoy, and Joel~A Tropp.
\newblock Living on the edge: Phase transitions in convex programs with random data.
\newblock {\em Information and Inference: A Journal of the IMA}, 3(3):224--294, 2014.

\bibitem[Bar16]{Bar}
Alexander Barvinok.
\newblock {\em Combinatorics and complexity of partition functions}, volume~30 of {\em Algorithms and Combinatorics}.
\newblock Springer, Cham, 2016.

\bibitem[BGL14]{BGL-book}
Dominique Bakry, Ivan Gentil, and Michel Ledoux.
\newblock {\em Analysis and geometry of {M}arkov diffusion operators}, volume 348 of {\em Grundlehren der mathematischen Wissenschaften [Fundamental Principles of Mathematical Sciences]}.
\newblock Springer, Cham, 2014.

\bibitem[BL14]{BL-Ber}
Witold Bednorz and Rafa{\l} Lata{\l}a.
\newblock On the boundedness of {B}ernoulli processes.
\newblock {\em Ann. of Math. (2)}, 180(3):1167--1203, 2014.

\bibitem[BM02]{bartlett2002rademacher}
Peter~L Bartlett and Shahar Mendelson.
\newblock Rademacher and gaussian complexities: Risk bounds and structural results.
\newblock {\em Journal of Machine Learning Research}, 3(Nov):463--482, 2002.

\bibitem[CGS11]{CGS-book}
Louis H.~Y. Chen, Larry Goldstein, and Qi-Man Shao.
\newblock {\em Normal approximation by {S}tein's method}.
\newblock Probability and its Applications (New York). Springer, Heidelberg, 2011.

\bibitem[CH06]{CH-univ}
Philippe Carmona and Yueyun Hu.
\newblock Universality in {S}herrington--{K}irkpatrick's {S}pin {G}lass {M}odel.
\newblock In {\em Annales de l'Institut Henri Poincare (B) Probability and Statistics}, volume~42, pages 215--222. Elsevier, 2006.

\bibitem[Cha05]{chat-error}
Sourav Chatterjee.
\newblock An {E}rror {B}ound in the {S}udakov-{F}ernique {I}nequality.
\newblock {\em arXiv preprint math/0510424}, 2005.

\bibitem[Cha14]{Cha-icm}
Sourav Chatterjee.
\newblock A short survey of {S}tein's method.
\newblock In {\em Proceedings of the {I}nternational {C}ongress of {M}athematicians---{S}eoul 2014. {V}ol. {IV}}, pages 1--24. Kyung Moon Sa, Seoul, 2014.

\bibitem[Kla06]{K}
Bo'az Klartag.
\newblock On convex perturbations with a bounded isotropic constant.
\newblock {\em Geom. Funct. Anal.}, 16(6):1274--1290, 2006.

\bibitem[KM11]{korada2011applications}
Satish~Babu Korada and Andrea Montanari.
\newblock Applications of the {L}indeberg principle in communications and statistical learning.
\newblock {\em IEEE transactions on information theory}, 57(4):2440--2450, 2011.

\bibitem[KM12]{KM}
Bo'az Klartag and Emanuel Milman.
\newblock Centroid bodies and the logarithmic {L}aplace transform---a unified approach.
\newblock {\em J. Funct. Anal.}, 262(1):10--34, 2012.

\bibitem[LW08]{LW}
Rafa{\l} Lata{\l}a and Jakub~Onufry Wojtaszczyk.
\newblock On the infimum convolution inequality.
\newblock {\em Studia Math.}, 189(2):147--187, 2008.

\bibitem[NP12]{nourdin2012normal}
Ivan Nourdin and Giovanni Peccati.
\newblock {\em Normal approximations with Malliavin calculus: from Stein's method to universality}, volume 192.
\newblock Cambridge University Press, 2012.

\bibitem[OT18]{oymak2018universality}
Samet Oymak and Joel~A Tropp.
\newblock Universality laws for randomized dimension reduction, with applications.
\newblock {\em Information and Inference: A Journal of the IMA}, 7(3):337--446, 2018.

\bibitem[RV08]{RV-recon}
Mark Rudelson and Roman Vershynin.
\newblock On sparse reconstruction from {F}ourier and {G}aussian measurements.
\newblock {\em Comm. Pure Appl. Math.}, 61(8):1025--1045, 2008.

\bibitem[Ste72]{Stein1972ABF}
Charles Stein.
\newblock A bound for the error in the normal approximation to the distribution of a sum of dependent random variables.
\newblock In {\em Proceedings of the {S}ixth {B}erkeley {S}ymposium on {M}athematical {S}tatistics and {P}robability ({U}niv. {C}alifornia, {B}erkeley, {C}alif., 1970/1971), {V}ol. {II}: {P}robability theory}, pages 583--602. Univ. California Press, Berkeley, CA, 1972.

\bibitem[Tal87]{Tal-reg}
Michel Talagrand.
\newblock Regularity of {G}aussian processes.
\newblock {\em Acta Math.}, 159(1-2):99--149, 1987.

\bibitem[Tal96]{talagrand1996majorizing}
Michel Talagrand.
\newblock Majorizing measures: the generic chaining.
\newblock {\em The Annals of Probability}, 24(3):1049--1103, 1996.

\bibitem[Tal01]{talagrand2001majorizing}
Michel Talagrand.
\newblock Majorizing measures without measures.
\newblock {\em Annals of probability}, pages 411--417, 2001.

\bibitem[Tal02]{Tal-Gibbs}
Michel Talagrand.
\newblock Gaussian averages, {B}ernoulli averages, and {G}ibbs' measures.
\newblock {\em Random Structures \& Algorithms}, 21(3-4):197--204, 2002.

\bibitem[Tal03]{Tal-spin}
Michel Talagrand.
\newblock {\em Spin glasses: a challenge for mathematicians}, volume~46 of {\em Ergebnisse der Mathematik und ihrer Grenzgebiete. 3. Folge. A Series of Modern Surveys in Mathematics. Cavity and mean field models. [Results in Mathematics and Related Areas. 3rd Series. A Series of Modern Surveys in Mathematics]}.
\newblock Springer-Verlag, Berlin, 2003.

\bibitem[Tal21]{Tal-book}
Michel Talagrand.
\newblock {\em Upper and lower bounds for stochastic processes---decomposition theorems}, volume~60 of {\em Ergebnisse der Mathematik und ihrer Grenzgebiete. 3. Folge. A Series of Modern Surveys in Mathematics [Results in Mathematics and Related Areas. 3rd Series. A Series of Modern Surveys in Mathematics]}.
\newblock Springer, Cham, second edition, \copyright 2021.

\bibitem[TH15]{7282516}
Christos Thrampoulidis and Babak Hassibi.
\newblock Isotropically random orthogonal matrices: Performance of lasso and minimum conic singular values.
\newblock In {\em 2015 IEEE International Symposium on Information Theory (ISIT)}, pages 556--560, 2015.

\bibitem[Ver18]{vershynin2018high}
Roman Vershynin.
\newblock {\em High-dimensional probability: An introduction with applications in data science}, volume~47.
\newblock Cambridge university press, 2018.

\end{thebibliography}

\end{document}